\documentclass[twoside,a4paper,12pt,centertags]{amsart}

\usepackage{amsmath,amssymb,verbatim,vmargin}
\usepackage[all]{xy}
\usepackage[bookmarks=true]{hyperref}   % Hyperref in DVI and PDF

\theoremstyle{plain}
\newtheorem{thm}{Theorem}[section]
\newtheorem{lem}[thm]{Lemma}
\newtheorem{cor}[thm]{Corollary}
\newtheorem{prop}[thm]{Proposition}

\theoremstyle{definition}
\newtheorem{defn}[thm]{Definition}
\newtheorem{rem}[thm]{Remark}

\newtheorem{ex}[thm]{Example}

\title[Hilbert transforms and the Cauchy integral in euclidean space]
{Hilbert transforms and the Cauchy integral in euclidean space}
\author{Andreas Axelsson} \author{Kit Ian Kou} \author{Tao Qian}
\address{Andreas Axelsson, Matematiska institutionen, Stockholms universitet, 106 91 Stockholm, Sweden}
\email{andax@math.su.se}
\address{Kit Ian Kou and Tao Qian, Department of Mathematics, University of Macau, Taipa, Macau, China}
\email{kikou@umac.mo, fsttq@umac.mo}

\newcommand{\dirac}{{\mathbf D}}
\newcommand{\del}{\delta}
\newcommand{\wedg}{\mathbin{\scriptstyle{\wedge}}}

\newcommand{\lctr}{\mathbin{\lrcorner}}
\newcommand{\rctr}{\mathbin{\llcorner}}
\newcommand{\clp}{\mathbin{\scriptscriptstyle{\triangle}}}

\newcommand{\scl}[2]{\langle #1,#2 \rangle}

\newcommand{\inv}[1]{{\widehat{#1}}}

\mathchardef\semic="303B

\newcommand{\R}{{\mathbf R}}
\newcommand{\C}{{\mathbf C}}

\newcommand{\mX}{{\mathcal X}}
\newcommand{\mY}{{\mathcal Y}}

\DeclareMathOperator{\re}{Re}
\newcommand{\im}{\text{{\rm Im}}\,}
\newcommand{\sett}[2]{ \{ #1 \, \semic \, #2 \} }

\newcommand{\supp}{\text{{\rm supp}}\,}
\newcommand{\dist}{\text{{\rm dist}}\,}

\newcommand{\nul}{\textsf{N}}
\newcommand{\ran}{\textsf{R}}
\newcommand{\dom}{\textsf{D}}

\newcommand{\clos}[1]{\overline{#1}}

\newcommand{\pd}{\partial}

\newcommand{\barint}{\mbox{$ave \int$}}
\def\barint_#1{\mathchoice
            {\mathop{\vrule width 6pt
height 3 pt depth -2.5pt
                    \kern -8.8pt
\intop}\nolimits_{#1}}%
            {\mathop{\vrule width 5pt height
3 pt depth -2.6pt
                    \kern -6.5pt
\intop}\nolimits_{#1}}%
            {\mathop{\vrule width 5pt height
3 pt depth -2.6pt
                    \kern -6pt
\intop}\nolimits_{#1}}%
            {\mathop{\vrule width 5pt height
3 pt depth -2.6pt
          \kern -6pt \intop}\nolimits_{#1}}}

\newcommand{\bdy}{\partial}
\newcommand{\ud}{\underline{d}}
\newcommand{\udel}{\underline{\delta}}
\newcommand{\pv}{\text{p.v.\!}}
\newcommand{\nindex}{{\overline n}}

\begin{document}

\maketitle

\begin{abstract}
   We generalize the notions of harmonic conjugate functions and Hilbert transforms
   to higher dimensional euclidean spaces, in the setting of differential forms and 
   the Hodge-Dirac system. 
   These harmonic conjugates are in general far from being unique, but under suitable boundary
   conditions we prove existence and uniqueness of conjugates.
   The proof also yields invertibility results for a new class of generalized double layer potential operators 
   on Lipschitz surfaces and boundedness of related Hilbert transforms.
\end{abstract}



\section{Introduction}

This paper considers higher dimensional analogues of the concept of harmonic conjugate functions in the plane.
We first review the situation for plane domains.
Let $D$ be a simply connected domain in $\R^2=\C$. Then, given a harmonic function $u(z)$ in $D$,
a {\em harmonic conjugate} to $u(z)$ is a second harmonic function $v(z)$ in $D$ such that
$f(z)=u(z)+i v(z)$ is an analytic function, i.e. satisfies Cauchy--Riemann's equations.
The function $v(z)$ exists and is unique modulo constants.
For example, if one requires that $v$ vanishes at some fixed point in $D$, then we get a well 
defined map $u\mapsto v$.
Since harmonic functions are in one-to-one correspondence with their boundary values, 
this defines the {\em Hilbert transform} 
$$
  H_D: u|_{\bdy D}\mapsto v|_{\bdy D}
$$ 
for the domain.
The Hilbert transform for a domain $D\subset \R^2$ concerns only functions in $D$, mapping the real part
of an analytic function to its imaginary part.
In contrast, the Cauchy integral concerns the relation between analytic functions in $D$
and $\C\setminus \clos D$.
This is best explained through Hardy spaces/projections. Given a function $h:\bdy D\rightarrow \C$,
we form the Cauchy integral
$$
  F(z)= C^\pm h(z):=\pm \frac 1{2\pi i}\int_{\bdy D}\frac {h(w)}{w-z} dw,\qquad z\in D^\pm,
$$
where $D^+:= D$ and $D^-:=\C\setminus\clos D$.
Taking traces, one obtains two boundary functions $f^\pm(\zeta):= \lim_{z\rightarrow \zeta, z\in D^\pm} F(z)$ such that
$f^++ f^-=h$ on $\bdy D$. 
The Cauchy integral acts by projection onto the two complementary 
Hardy subspaces of the boundary function space,
consisting of traces $f^+$ and $f^-$ respectively.
Discarding the exterior Hardy function $f^-$ and considering $f:= f^+= u+iv$, well posedness of the 
classical Hilbert boundary value problem (BVP) for analytic functions shows that $f$ is
in one-to-one correspondence with its real part $u$, as well as with its imaginary part $v$.
(In this introduction, we neglect technical details like regularity assumptions
on $\bdy D$ and $h$, as well as the fact that maps normally are Fredholm, not exact 
isomorphisms, in order to explain the main ideas.)
From the Hilbert BVP it can be shown that if $h$ is required to be real valued, then $h$ is also 
in one-to-one correspondence with the trace of its interior Cauchy integral 
$f=C^+ h|_{\bdy D}= \tfrac 12h+ \tfrac 1{2\pi i}\pv \int_{\bdy D}h(w)/(w-z)dw$.
Thus in the diagram
$$
\xymatrix{ 
  & f \ar[r]^{\im} \ar[rd]^{\re}& v \\
  h \ar[ru]^{C^+|_{\bdy D}} \ar[rr]_T & & u
 }
$$
the {\em double layer potential operator}
$$
  Th(z):= \re( C^+h|_{\bdy D}(z) )= \frac 12 h(z)+ \frac 1{2\pi}\pv \int_{\bdy D}\im\Big(\frac {dw}{w-z}\Big)h(w),
  \qquad z\in \bdy D,
$$
is an isomorphism. This allows calculation of the Hilbert transform
as
\begin{equation}  \label{eq:hilbertfactor}
  v= H_Du= \im( C^+( T^{-1}u )|_{\bdy D}).
\end{equation}
Note that even though the Cauchy integral, for all $D$, uses the restriction of $1/(w-z)$ as kernel, the 
kernel of the Hilbert transform depends heavily on $D$ due to the factor $T^{-1}$.
Unfortunately, in the literature the Hilbert transform is often incorrectly identified with the Cauchy integral, since it  happens to coincide with the (imaginary part of) the Cauchy integral for half planes and disks. 
However, we emphasize that for
all other domains, these two operators are not the same.

The aim of this paper is to show the existence and $L_2$-boundedness of higher dimensional
Hilbert transforms for Lipschitz domains in $\R^n$. 
These Hilbert transforms are derived from a Cauchy integral as in (\ref{eq:hilbertfactor}).
In the plane, the double layer potential operator $T$ is the {\em compression} of the Cauchy integral
$C^+|_{\bdy D}$ to the subspace of real valued functions, i.e. $T$ is $C^+|_{\bdy D}$ restricted
to the subspace of real functions, followed by projection back onto this subspace.
In higher dimensions, there is a canonical Cauchy integral as well, but there are various 
natural subspaces to compress it to, to obtain generalizations of the double layer potential operator.
It is of importance to establish boundedness and invertibility of such operators, as seen
for example in (\ref{eq:hilbertfactor}).
In connection with boundary value problems for Dirac operators, see the works \cite{Ax1, Ax2, Ax} by the first author
and Remark~\ref{rem:compressions} for more details,
one type of compressed Cauchy integrals plays a central role.
The compressed Cauchy integrals used in this paper to calculate higher dimensional Hilbert transforms
have not been studied before, to the authors knowledge.

We next explain the higher dimensional concepts of harmonic conjugate functions and Hilbert
transforms that we consider in this paper.
The plane domain $D$ is replaced by a domain in $\R^n$ and the Cauchy--Riemann system is
generalized to the {\em Hodge--Dirac system} $(d+\del) F(x)=0$.
Functions $F:D\rightarrow \wedge\R^n$ now in general take values in the full $2^n$-dimensional exterior algebra
$$
  \wedge \R^n= \wedge^0\R^n\oplus \wedge^1\R^n\oplus \wedge^2\R^n\oplus\ldots\oplus  
  \wedge^{n-1}\R^n\oplus \wedge^n\R^n,
$$
and $d$ and $\del$ denote the exterior and interior derivative operators
\begin{align*}
  dF(x) &= \nabla\wedg F(x)= \sum_{j=1}^n e_j\wedg \pd_j F(x), \\
  \del F(x) &= \nabla\lctr F(x)= \sum_{j=1}^n e_j\lctr \pd_j F(x),
\end{align*}
where the exterior and interior products are dual in the sense that $\scl x{u\wedg y}=\scl{u\lctr x}y$ for all 
multivectors $u,x,y\in \wedge \R^n$.
Multivector fields $F$ taking values in the subspace $\wedge^k \R^n$ we refer to as $k$-vector fields.
The Hodge--Dirac equation entails a coupling between the different $k$-vector parts of $F$.
If $F=\sum_k F_k$, where $F_k:D\rightarrow \wedge^k\R^n$, the differential operators map
\begin{equation}  \label{eq:complex}
\xymatrix{ 
  \cdots & \ar@<-1ex>[l];[]_>>>d  \ar@<-1ex>[l]_>>>\del \wedge^{k-2}\R^n \ar@<1ex>[r]^>>>d & \ar@<1ex>[l]^>>>\del \wedge^{k-1}\R^n \ar@<-1ex>[r]_>>>d & \ar@<-1ex>[l]_>>>\del \wedge^{k}\R^n \ar@<1ex>[r]^>>>d & \ar@<1ex>[l]^>>>\del \wedge^{k+1}\R^n \ar@<-1ex>[r]_>>>d & \ar@<-1ex>[l]_>>>\del \wedge^{k+2}\R^n \ar@<1ex>[r]^>>>d & \ar@<1ex>[l]^>>>\del \cdots
 }
\end{equation}
so that the Hodge--Dirac equation is equivalent to the system of equations $dF_{k-1}= -\del F_{k+1}$, $0\le k\le n$.
If $F$ is {\em monogenic}, i.e. satisfies the Hodge--Dirac equation, then it is {\em harmonic}, i.e. satisfies 
$\Delta F= (d+\del)^2 F= (\del d+d\del) F=0$,
or equivalently each of $F$'s $2^n$ scalar component functions is harmonic. 
Recall that $d^2=\del^2=0$.

Fix $0\le k\le n$ and consider a harmonic $k$-vector field $U: D\rightarrow \wedge^k\R^n$. 
We say that $V_1: D\rightarrow \wedge^{k-2}\R^n$ and $V_2: D\rightarrow \wedge^{k+2}\R^n$ form a pair
of {\em harmonic conjugates} to $U$ if $(d+\del)(V_1+U+ V_2)=0$, or
equivalently $d U=-\del V_2$, $\del U=-d V_1$ and $d V_2=\del V_1=0$.
\begin{ex}
When $n=2$ and $k=0$ this reduces to the classical situation.
Indeed, consider an analytic function $f=u+iv$.
We identify $\R\approx\wedge^0\R^2$ and $i\R\approx\wedge^2\R^2$.
Fixing an ON-basis $\{e_1, e_2\}$ for $\R^2$, we identify $i\approx e_1\wedg e_2$.
Then
\begin{multline*}
  (d+\del)(u+iv)=(e_1\wedg \pd_1 u+ e_2\wedg\pd_2 u)+(e_1\lctr i\pd_1 v + e_2\lctr i\pd_2 v) \\
  = (\pd_1 u-\pd_2 v)e_1+ (\pd_2 u+ \pd_1 v)e_2=0
\end{multline*}
coincides with the Cauchy--Riemann equations.
Hence $v$ is a classical harmonic conjugate to $u$ if and only if $V_2=v\, e_1\wedg e_2: D\to \wedge^2\R^2$
is a harmonic conjugate to $u$ in the sense of the Hodge--Dirac system.
(In this case, the harmonic conjugate $V_1$ vanishes.)
\end{ex}
In the general case, we observe that a necessary condition for such $V_1$, $V_2$ to exist is that 
$U$ is {\em two-sided harmonic}, i.e. $\del d U=0= d\del U$. We also observe that $V_1$, $V_2$ are only 
well defined modulo {\em two-sided monogenic fields}, i.e. the differences $V_1-V_1'$ and $V_2-V_2'$
of two sets of harmonic conjugates satisfy $d(V_i-V_i')=0= \del(V_i-V_i')$, $i=1,2$.
When $1\le k \le n-1$, the two-sided monogenic $k$-vector fields form an infinite dimensional space
(see Corollary~\ref{cor:infdim}).
Thus, in order to obtain a uniquely defined higher dimensional Hilbert transform, further conditions need 
to be imposed on $V_1$ and $V_2$, so that there is a well defined map
$$
  U\longmapsto V_1,V_2.
$$
In this paper we consider one possible such further condition on harmonic conjugate functions,
which extends the above technique of calculating conjugates
with the Cauchy integral and double layer potential operators to higher dimension.
Under this further condition we say that
the harmonic conjugate functions are of {\em Cauchy type}.
Since all component functions of $U,V_1$ and $V_2$ are harmonic, these fields are in one-to-one 
correspondence with their trace on $\partial D$. Thus, equivalently we will have a 
Hilbert type transform $U|_{\partial D}\mapsto V_1|_{\partial D}, V_2|_{\partial D}$ for $D$.

The outline of the paper is as follows. In Section~\ref{sec:scalarHilbert} we 
introduce the higher dimensional Cauchy integral associated with the Hodge--Dirac system,
and prove existence and uniqueness results for Cauchy type conjugates to scalar functions,
i.e. $k=0$ or $k=n$.
This amounts to proving invertibility of the classical double layer potential
operator for domains in $\R^n$. 
Both the boundedness and invertibility of this singular integral operator on $L_p$-spaces on Lipschitz 
boundaries ($1<p<\infty$ for boundedness and $2-\epsilon<p<\infty$ for invertibility)
are deep results, but are by now well known facts.

Section~\ref{sec:kHilbert} contains the main new results of the paper, Theorem~\ref{thm:maink}, and establishes 
existence and uniqueness of Cauchy type conjugates to $k$-vector fields,
$1\le k\le n-1$.
This general case is more involved than the scalar case, since the generalized 
double layer potential operators which appear will not in general 
be invertible, not even Fredholm, as they have infinite dimensional null spaces and cokernels. 
However, using the theory of boundary value problems for Dirac operators 
(which is reviewed in Section~\ref{sec:bvps}), we manage
to show invertibility of the operator acting from a complement of the null space to its range,
in a natural $L_2$-based Hilbert space.

In the final Section~\ref{sec:hodge} we illustrate the non-uniqueness of harmonic conjugate functions in higher
dimensional euclidean spaces by constructing different conjugate functions which are 
not in general the Cauchy type conjugates.
This second construction is based on the theory of Hodge decompositions, and 
the obtained harmonic conjugate functions are said to be of {\em Hodge type}. 

In the literature, various generalizations of harmonic conjugate functions 
to higher dimensional euclidean spaces can be found. 
A classical generalization for the upper half space, using divergence and curl free vector fields,
was introduced in harmonic analysis by Stein and Weiss~\cite{SW}, see 
Stein~\cite{S} and Example~\ref{ex:cauchyconj}(1).
A generalization more similar to our construction is due to Ar\u zanyh~\cite{Arz},
who studied two-forms $B:D\rightarrow \wedge^2\R^3$ conjugate to scalar functions in three
dimensional space.
In the setting of Clifford analysis, without dealing directly with the more fundamental differential operators
$d$ and $\del$, there is work on Hilbert transforms and harmonic conjugate functions in euclidean space 
by Brackx, De Knock, De Schepper and Eelbode \cite{BKSE}.
See also the references therein for calculations on special domains like the unit ball.

\section{Hilbert transforms for scalar functions}   \label{sec:scalarHilbert}

Writing $\dirac := d+\del$ for the Hodge--Dirac operator, this is an elliptic first order partial
differential operator whose square is the Hodge--Laplace operator $\dirac^2=\Delta$.
Just like the exterior and interior differentiation use the exterior and interior products,
the Hodge--Dirac operator uses the Clifford product as 
$$
    \dirac F(x) = \nabla\clp F(x)= \sum_{j=1}^n e_j\clp \pd_j F(x).
$$
The Clifford product is the unique associative algebra product $\clp$ on $\wedge\R^n$, 
with identity $1\in\wedge^0\R^n=\R$, such that
\begin{equation}   \label{eq:cliffdefn}
   v\clp w= v\lctr w+ v\wedg w,\qquad w\clp v= w\rctr v+ w\wedg v
\end{equation}
for all vectors $v\in \wedge^1\R^n=\R^n$ and all multivectors $w\in \wedge\R^n$.
The main difference between the complex product in $\R^2$ and its higher dimensional
analogue in $\R^n$, the Clifford product, is that the latter is non-commutative. 
Here $\lctr$ and $\rctr$ denote the left and right interior products, defined as the operations
adjoint to left and right exterior multiplication, i.e.
\begin{equation}    \label{eq:intprod}
  \scl{w\lctr x}{y}=\scl x{w\wedg y},\qquad \scl{x\rctr w}{y}= \scl x{y\wedg w},\qquad w,x,y\in\wedge \R^n.
\end{equation}
Following standard notation, we write $w_1\clp w_2=: w_1w_2$ for short.
Important to this paper is the following mapping property of the Clifford product. If $v\in\wedge^1\R^n=\R^n$
and $w\in\wedge^k\R^n$, then
$$
  vw \in\wedge^{k-1}\R^n\oplus \wedge^{k+1}\R^n.
$$
This is clear from (\ref{eq:cliffdefn}).

Concretely, if $\{e_i\}_{i=1}^n$ denotes the standard ON-basis for $\R^n$, the induced ON-basis
for the space of $k$-vectors $\wedge^k\R^n$ is $\{e_s\}_{|s|=k}$, and in total
$\{e_s\}_{s\subset\nindex}$ is the induced basis for $\wedge\R^n$, where
$\nindex:= \{1,2,\ldots, n\}$.
For a subset $s=\{s_1,\ldots, s_k\}\subset\nindex$, where $1\le s_1<\ldots <s_k\le n$, we write
$e_s:= e_{s_1}\wedg\ldots\wedg e_{s_k}$.
In this induced basis the exterior, interior and Clifford products are
\begin{align*}
  e_s\wedg e_t  = \epsilon(s,t)\, e_{s\cup t},&\qquad \text{if } s\cap t=\emptyset, \text{ and otherwise } 0, \\
  e_s\lctr e_t  = \epsilon(s,t\setminus s) \, e_{t\setminus s},&\qquad \text{if } s\subset t, \text{ and otherwise } 0, \\
  e_s\rctr e_t  = \epsilon(s\setminus t,t) \, e_{s\setminus t},&\qquad \text{if } t\subset s, \text{ and otherwise } 0, \\
  e_se_t  = \epsilon(s,t)\, e_{t\clp s},&
\end{align*}
where the permutation sign is $\epsilon(s,t):= (-1)^{|\sett{(s_i,t_j)\in s\times t}{s_i>t_j}|}=\pm 1$
and $s\clp t:= (s\setminus t)\cup (t\setminus s)$ denotes the symmetric difference of index sets.

The radial vector field 
$$
  E(x):= \frac 1{\sigma_{n-1}}\frac x{|x|^n}, \qquad x\in \R^n,
$$
where $\sigma_{n-1}$ is the area of the unit sphere $S^{n-1}$,
is divergence and curl free for $x\ne 0$, and is a fundamental solution to the Hodge--Dirac operator.
Using the associativity of the Clifford product, Stokes' theorem gives a Cauchy type integral formula
\begin{equation}    \label{eq:cauchyintegral}
  F(x)=\int_{\bdy D} E(y-x) \nu(y) F(y) d\sigma(y), \qquad x\in D,
\end{equation}
for {\em monogenic fields} $F:\clos D\rightarrow \wedge\R^n$, i.e. solutions to $\dirac F=0$.
Here $\nu(y)$ denotes the unit normal vector field on $\bdy D$ pointing outward from $D$,
and $d\sigma$ denotes scalar surface measure on $\bdy D$.
Note that $\nu(y)$ must be placed between the two factors because of non-commutativity.

\begin{ex}
  If $n=2$ and $F:\R^2\approx\C\rightarrow \C\approx\wedge^0\R^2\oplus \wedge^2\R^2$, 
  the Cauchy formula (\ref{eq:cauchyintegral}) reduces to the standard one. 
  Indeed, 
$$
  E(y-x)\nu(y)d\sigma(y)= \frac 1{2\pi i}\frac{dw}{w-z}
$$
if we identify $i=e_1e_2$, $w= e_1 y$, $z= e_1 x$ and $dw/i= e_1 \nu d\sigma$.
\end{ex}

With these algebraic preliminaries, we next turn to analysis.
For estimates, we use the notation $X\approx Y$ to mean that there exists a constant $C$, independent
of the variables in the estimate, such that $C^{-1}Y\le X\le CY$. Similarly $X\lesssim Y$ means that $X\le CY$.

To avoid topological technicalities, we shall restrict attention to the following two
types of strongly Lipschitz domains $D^\pm \subset \R^n$ (i.e. domains whose boundaries are locally graphs of Lipschitz functions).
We write $\Sigma:= \bdy D^+=\bdy D^-$ for the Lipschitz boundary between $D^+$ and 
$D^-= \R^n\setminus\clos D^+$. By $D$ we denote either $D^+$ or $D^-$.
\begin{itemize}
\item A {\em graph domain} $D^+= \sett{x}{x_n>\phi(x_1,\ldots, n_{n-1})}$ above the graph of a Lipschitz 
regular function $\phi:\R^{n-1}\rightarrow \R$. Here $D^-$ denotes the domain below the graph.
\item An {\em interior domain} $D^+$, being a bounded domain which is Lipschitz diffeomorphic to the unit
ball, and whose boundary is locally the graph of a Lipschitz function (in suitable ON-bases).
The {\em exterior domain} $D^-$ is the interior of the unbounded complement of $D^+$.
\end{itemize}
The unit normal vector field $\nu(y)$ on $\Sigma$ is always assumed to point into $D^-$, i.e.
the region below the graph or into the exterior domain.
We define non-tangential approach regions $\gamma(y)\subset D$, $y\in\Sigma$, 
for these Lipschitz boundaries.
For graph domains $D^\pm$, fix $c_1$ greater than the Lipschitz constant for $\Sigma$, and let
$$
  \gamma(y)=\gamma(y,D^\pm):= \sett{(x', x_n)\in\R^{n-1}\times \R}{\pm(x_n-y_n)>c_1|x'-y'|}, 
$$
for $y=(y', y_n)\in \Sigma$.
For exterior and interior domains, and $y\in \Sigma$, consider the coordinate system 
around $y$ in a neighbourhood of which $\Sigma$ is a Lipschitz graph.
The approach region $\gamma(y,D^\pm)$ is defined as the truncated part of the cone, where
$\dist(x,y)<c_2$ and $c_2$ denotes a sufficiently small constant.

The boundary function spaces we use are the spaces $L_p(\Sigma;\wedge)$, where 
$\wedge:= \wedge\R^n$.
For a field $F$ in $D^\pm$, define its {\em non-tangential maximal function} 
$$
  N_*(F)(y):= \sup_{x\in \gamma(y,D^\pm)}|F(x)|,\qquad y\in \Sigma.
$$
A fundamental theorem in harmonic analysis and singular integral theory due to
Coifman, McIntosh and Meyer~\cite{CMcM}, states that the Cauchy integral is bounded
on $L_p(\Sigma;\wedge)$ on any Lipschitz surface $\Sigma$.
By surface, we shall mean a hypersurface in $\R^n$.

\begin{thm}    \label{thm:cmcm}
  Let $D^\pm$ be Lipschitz graph, interior or exterior domains, and fix $1<p<\infty$. 
  Let $h\in L_p(\Sigma;\wedge)$ and define the monogenic field
$$
  C^\pm h(x):= \pm \int_\Sigma E(y-x)\nu(y) h(y)d\sigma(y),\qquad x\in D^\pm.
$$
Then $\|N_*(C^+ h)\|_p+ \|N_*(C^-h)\|_p\le C\|h\|_p$ for some $C<\infty$ depending only on $p$ and 
the Lipschitz constants for the graphs describing $\Sigma$. 

The principal value Cauchy integral
$$
  Eh(x):= 2\pv \int_\Sigma E(y-x)\nu(y) h(y) d\sigma(y),\qquad x\in \Sigma,
$$
exists a.e. and defines a bounded operator $E:L_p(\Sigma;\wedge)\rightarrow L_p(\Sigma;\wedge)$
such that $E^2=I$.
The boundary traces $f^+(z):= \lim_{x\rightarrow z, x\in \gamma(z, D^+)} C^+h(x)$
and $f^-(z):= \lim_{x\rightarrow z, x\in \gamma(z, D^-)} C^-h(x)$ exist for a.a. $z\in\Sigma$ and in $L_p$, and
$$
  E^+h:= f^+=\tfrac 12(h+Eh) \qquad \text{and}\qquad E^-h:= f^-=\tfrac 12(h-Eh)
$$
define $L_p$-bounded projection operators. 
\end{thm}

Let $h\in L_p(\Sigma;\wedge^k)$ be a $k$-vector field and consider harmonic
conjugates in $D^+$. 
The Cauchy integral produces a monogenic field 
$$
  F=V_1+U+V_2= C^+h: D^+\longrightarrow \wedge^{k-2}\R^n\oplus\wedge^k\R^n\oplus
  \wedge^{k+2}\R^n
$$
in $D^+$.
Indeed, the mapping properties of the Clifford product show that multiplication
with the normal vector gives $\nu h:\Sigma\rightarrow \wedge^{k-1}\R^n\oplus \wedge^{k+1}\R^n$,
and a similar splitting when multiplying with the vector $E(y-x)$ shows that
\begin{multline}   \label{eq:Cauchymapping}
  C^+h(x)= \int_\Sigma E(y-x)\lctr(\nu(y)\lctr h(y))\, d\sigma(y) \\
  +\int_\Sigma \Big( E(y-x)\wedg(\nu(y)\lctr h(y))+ E(y-x)\lctr(\nu(y)\wedg h(y)) \Big) \,d\sigma(y) \\
  +\int_\Sigma E(y-x)\wedg \nu(y) \wedg h(y) \, d\sigma(y) = V_1+U+V_2.
\end{multline}

\begin{defn}
  Given a two-sided harmonic field $U:D^\pm\rightarrow \wedge^k\R^n$, i.e. $\del dU=0= d\del U$,
  we say that $V_1,V_2$ are {\em Cauchy type harmonic conjugates} to $U$ if
  there exists $h:\Sigma\rightarrow\wedge^k\R^n$ such that $U= (C^\pm h)_k$, $V_1= (C^\pm h)_{k-2}$
  and $V_2=(C^\pm h)_{k+2}$, where subscript $k$ denotes the $k$-vector part of a multivector.
  We call $h$ the {\em dipole density} of the system $V_1, U, V_2$ of harmonic conjugate functions. 
  The corresponding map of boundary values
  $$
    U|_\Sigma\longmapsto V_1|_{\Sigma}, V_2|_{\Sigma},
  $$
  we refer to as the (Cauchy type) {\em Hilbert transform} for the domain $D$.
\end{defn}

The following theorem on existence and uniqueness of Cauchy type harmonic conjugates to 
scalar functions ($k=0,n$) is the main result of this section.

\begin{thm}   \label{thm:main}
  Let $D\subset \R^n$ be a Lipschitz graph, interior or exterior domain and assume that $2\le p<\infty$.
\begin{itemize}
\item[{\rm (i)}]
  Let $U:D \rightarrow \R=\wedge^0\R^n$ be a harmonic function such that $N_*(U)\in L_p(\Sigma)$.
  If $D$ is an exterior domain, also assume that $\lim_{x\rightarrow \infty} U=0$ and 
  has trace $u=U|_\Sigma$ such that $\int u\psi\, d\sigma(y) =0$, where $\psi$ is the function from 
  Theorem~\ref{thm:dahlbergkenig}.
  Then there is a unique Cauchy type harmonic conjugate $V=V_2:D\rightarrow \wedge^2\R^n$
  to $U$, and a dipole density $h\in L_p(\Sigma)$, such that 
$$
  \|N_*(V)\|_p+ \|h\|_p\lesssim \|N_*(U)\|_p.
$$
If $D$ is a graph or an interior domain, then $h$ is unique, and if $D$ is an exterior domain, then
$h$ is unique modulo constants.

\item[{\rm (ii)}]
  In the case $k=n$, (i) remains true when $U:D\rightarrow \wedge^0\R^n=\R$
  and $V=V_2:D\rightarrow \wedge^2\R^n$ are replaced by 
  $U:D \rightarrow \wedge^n\R^n\approx\R$ and $V=V_1:D\rightarrow \wedge^{n-2}\R^n$.
\end{itemize}
\end{thm}
The scalar cases $k=0$ and $k=n$ are significantly more straightforward than the non-scalar 
case $1\le k\le n-1$ (to be treated in Section~\ref{sec:kHilbert})
as they reduce to the question whether the classical double layer potential equations are invertible, as 
explained in the two dimensional case in the introduction.
On $\Sigma$, define the {\em principal value double layer potential operator}
\begin{equation}   \label{eq:dlpdefn}
  Kh(x):= 2\pv\int_\Sigma \scl{E(y-x)}{n(y)}h(y) \, d\sigma(y)= (E h(x))_0, 
\end{equation}
for $h:\Sigma\rightarrow\R, \, x\in \Sigma$.
The boundedness of $K$ in $L_p(\Sigma)$, $1<p<\infty$, is a direct consequence of 
Theorem~\ref{thm:cmcm}.
Invertibility of $I\pm K$ on the other hand, which the proof of Theorem~\ref{thm:main} uses, 
is not true for all $p$ on a general Lipschitz surface.
Invertibility in $L_2(\Sigma)$ was proved by Verchota~\cite{V} via the method of Rellich estimates.
Invertibility in the range $2<p<\infty$ was proved by Dahlberg and Kenig~\cite{DK} 
by atomic estimates in real Hardy space $H^1(\Sigma)$, duality and interpolation.

\begin{thm}   \label{thm:dahlbergkenig}
  Assume that $2\le p<\infty$ and let $\Sigma$ be a Lipschitz graph or the boundary of an interior / exterior domain.
  Then 
$$
  I+ K: L_p(\Sigma)\longrightarrow L_p(\Sigma): h\longmapsto  h+(Eh)_0=2(E^+ h)_0
$$ 
is an isomorphism. This is also true for 
$I-K: h\mapsto (E^-h)_0$ in the case of a graph domain.
In the case of an exterior domain, $I-K$ is a Fredholm operator with null space
consisting of constant functions and range consisting of all $u\in L_p(\Sigma)$
such that $\int_\Sigma u\psi=0$ for some $\psi\in L_q(\Sigma)$, where $1/p+1/q=1$.
\end{thm}
We remark that Theorems~\ref{thm:cmcm}, \ref{thm:main} and \ref{thm:dahlbergkenig} can be 
generalized to strongly Lipschitz domains
with more complicated topology. In this case existence and uniqueness of conjugates hold only modulo finite dimensional subspaces, 
and $I\pm K$ are Fredholm operators
with higher dimensional null spaces and cokernels.

\begin{proof}[Proof of Theorem~\ref{thm:main}]
   To prove (i), take $h\in L_p(\Sigma; \wedge^0)$ and define $\tilde U:= (C^\pm h)_0$ and 
   $\tilde u:= \tilde U|_\Sigma= \tfrac 12(I\pm K)h\in L_p(\Sigma; \wedge^0)$. 
  Theorem~\ref{thm:dahlbergkenig} determines uniquely $h$ such that $\tilde u=u$,
  possibly modulo constants in the case of an exterior domain.
  In any case,  this defines uniquely a Cauchy type harmonic conjugate $V:= (C^\pm h)_2$,
  since $C^-$ maps constants to zero in an exterior domain.
  
  To prove (ii), consider a system of Cauchy type harmonic conjugates
  $C^\pm h= V+ U$, where $h: \Sigma\rightarrow \wedge^n\R^n$,
  $V: D\rightarrow \wedge^{n-2}\R^{n}$ and $U: D\rightarrow \wedge^n\R^n$.
  Introduce the operator $U\mapsto Ue_\nindex$ (i.e. the Hodge star operator
  for differential forms up to a sign)
  which maps $k$-vectors to $n-k$-vectors.
  We have
  $C^\pm(he_\nindex)= Ue_\nindex+ Ve_\nindex$, where the functions
  take values in $\wedge^0\R^n$, $\wedge^0\R^n$ and $\wedge^2\R^n$ respectively.
  This reduces (ii) to (i), since $U\mapsto Ue_\nindex$ commutes with the Dirac operator
  and the Cauchy integral by the associativity of the Clifford product.
\end{proof} 

We note from the proof the following relations between Cauchy type harmonic conjugate
functions and Hilbert transforms of scalar functions, the Cauchy integral and 
the classical double layer potential operator.

\begin{cor}
  Let $V:D^+\rightarrow \wedge^2\R^n$ be the Cauchy type harmonic conjugate to
the harmonic function $U:D^+\rightarrow \wedge^0\R^n=\R$, with suitable estimates
of non-tangential maximal functions.
Then
$$
  U+V= 2C^+((I+K)^{-1}u),
$$
where $C^+$ is the (interior) Cauchy integral, $K$ is the double layer potential
operator, and $u=U|_\Sigma$.
Taking the trace $v=V|_\Sigma$ of the conjugate function, the Cauchy type Hilbert transform
of $u$ is
$$
  u\longmapsto v= (I+E)(I+K)^{-1}u.
$$
Replacing $C^+$, $I+K$ and $I+E$ with $C^-$, $I-K$ and $I-E$, the corresponding formulae hold
for the domain $D^-$.
\end{cor}

\begin{ex}   \label{ex:cauchyconj}
(1)
  If $D=\R^n_+$ is the upper half space with the flat boundary $\Sigma=\R^{n-1}$, 
  then $K=0$ since $E(y-x)$ in this case is orthogonal to $\nu(y)$.
  Thus $h=2u$ and 
$$
  V(x)= 2(C^+ u)_2(x)= e_n\wedg \frac 2{\sigma_{n-1}} \int_{\R^{n-1}} \frac{y-x}{|y-x|^n}u(y)\, dy, 
  \qquad x\in\R^n_+,
$$
 where $e_n$ is the basis vector normal to $\R^{n-1}$,
  so in this case the Hilbert transform $U|_{\R^{n-1}}\mapsto V|_{\R^{n-1}}$ for the upper half space 
  coincides with the bivector part of the principal value Cauchy integral $u\mapsto (Eu)_2$.

  A classical higher dimensional notion of harmonic conjugates, using divergence and curl free vector fields
  in the upper half space $\R^n_+$, is due to Stein and Weiss~\cite{SW}, and we refer to
  Stein's book \cite{S} for further details.
  The upper half space has the special property that the vector $e_n$ normal to $\bdy \R^n_+=\R^{n-1}$
  is constant. Split a vector field $F$ in $\R^n_+$ into normal and tangential parts as
$$
  F(x)= U(x)e_n + \tilde U (x),
$$
where $U$ is a scalar function.
Stein and Weiss consider the tangential vector field $\tilde U$ as a harmonic conjugate to $U$,
if $F$ is a divergence and curl free vector field, i.e. if $F$ is monogenic.
Since the Clifford product is associative, this is equivalent to
$U+ \tilde U e_n:\R^n_+\rightarrow \wedge^0\R^n\oplus \wedge^2\R^n$ being monogenic since
$$
  \dirac(U+\tilde U e_n)=\dirac((Ue_n+\tilde U)e_n)=(\dirac (Ue_n+\tilde U))e_n=0.
$$
Due to the very special geometry of $\R^n_+$, the bivector field $\tilde U e_n$ will in fact be the Cauchy type
harmonic conjugate to $U$. 
Indeed, if $V$ denotes the Cauchy type conjugate, then as noted in the introduction, the
difference $\tilde U e_n-V$ is a two-sided monogenic bivector field.
Moreover $V(x)= e_n\wedg \int_{\R^n} E(y-x)u(y)dy$, since the normal vector is constant.
Hence $e_n\wedg(\tilde U e_n-V)|_{\R^{n-1}}=0$.
From Theorem~\ref{thm:rellich} below we deduce that $\tilde U e_n-V=0$, since there are
no non-trivial monogenic field which are normal (in the sense of Definition~\ref{defn:normtang}) 
on the boundary.

Thus, from the above relation $\tilde U e_n= 2(C^+ u)_2$, the Stein--Weiss tangential vector field 
$\tilde U$ is seen to be
$$
  \tilde U(x)= \frac 2{\sigma_{n-1}} \int_{\R^n_+} \frac {x'-y}{|x-y|^n} u(y) dy,\qquad x=(x',x_n)\in\R^n_+,
$$
and taking the trace of $\tilde U$, 
the $n-1$ component functions of $\tilde U|_{\R^{n-1}}$ are the Riesz transforms of $u$.
  We remark that harmonic conjugates in the sense of Stein and Weiss do not generalize to more
  general domains $D$, since they depend on a canonical direction $e_n$.

(2)    
  If $D$ is the unit disk in the plane, then 
$$
  Kh(x)= \frac 1\pi\pv\int_\Sigma \frac {\scl{y-x}{y}}{|y-x|^2} h(y) \, d\sigma(y)= [h],
$$
where $[h]$ denotes the mean value of $h$, regarded as a constant function.
This gives $(I+K)^{-1}u= u-[u]/2$ and
\begin{multline*}
  V(x)= (2C^+(I+K)^{-1}u)_2(x)= (2C^+u-[u])_2= 2(C^+u)_2(x) \\ = 
  \frac 1\pi\int_\Sigma\frac {y\wedg x}{|y-x|^2} u(y)\, d\sigma(y)
  \longrightarrow (Eu)_2(e^{it})= i\frac 1{2\pi}\pv\int_0^{2\pi} \cot((t-s)/2) u(s)\, ds,
\end{multline*}
when $x\rightarrow e^{it}$, $y=e^{is}$, $i=e_1e_2$, and the limit is pointwise a.e. and in $L_p$, $1<p<\infty$.
Also here the Hilbert transform coincides with the imaginary part of the principal value Cauchy integral.
Note that this is far from being the case for more general domains, not even when $n=2$,
or for higher dimensional spheres as we shall see below.
For the exterior of the unit circle, we see that $(I-K)h=h-[h]$. Thus $I-K$ has constants as null space
and its range consists of all function with mean value zero. 
Hence the function $\psi$ in Theorem~\ref{thm:main}, orthogonal to the range, is a constant function.

(3)
For the unit sphere $\Sigma$ in $\R^n$, $K$ is the operator
$$
  Kh(x)=\frac 1{2^{n/2}\sigma_{n-1}}\int_\Sigma\frac{h(y)}{(1-\scl xy)^{n/2-1}}\,  d\sigma(y).
$$
This is weakly singular and therefore compact on all $L_p$ spaces. 
In fact, this is true whenever $\Sigma$ is a smooth surface, since the normal vector will be approximately orthogonal to $E(y-x)$ in the kernel when $y$ is close to $x$.
Through a limiting argument,
it was proved by Fabes, Jodeit and Rivi\`ere~\cite{FJR} that $K$ is a compact operator on 
$L_p$, $1<p<\infty$, whenever $\Sigma$ is a bounded $C^1$ regular surface. 
\end{ex}

\section{Dirac boundary value problems}     \label{sec:bvps}

In this section we describe parts of the operator theoretic framework for Dirac boundary 
value problems, developed by the first author in his PhD thesis \cite{Ax}, which we use
in Section~\ref{sec:kHilbert} to extend Theorem~\ref{thm:main}
to general $k$-vector fields.
In this section, we focus on explaining the main ideas of proofs, but give references to full proofs.
By $\nul(T)$, $\ran(T)$ and $\dom(T)$ we denote the null space, range and domain of an operator $T$.

The basic picture is that the boundary function space
\begin{equation}   \label{eq:homoparts}
  L_2(\Sigma):= L_2(\Sigma;\wedge)
  =L_2(\Sigma;\wedge^0)\oplus L_2(\Sigma;\wedge^1)\oplus\ldots\oplus
  L_2(\Sigma;\wedge^{n-1})\oplus L_2(\Sigma;\wedge^n)
\end{equation}
splits in two different ways into pairs of complementary closed subspaces
$$
  L_2(\Sigma)= E^+L_2\oplus E^-L_2= N^+L_2\oplus N^-L_2.
$$
In the first splitting $L_2=E^+L_2\oplus E^-L_2$, the subspaces $E^\pm L_2$ 
denote the Hardy type subspaces associated with the Dirac equation, i.e. $E^+L_2$ consists of traces of monogenic fields
in $D^+$ and $E^-L_2$ consists of traces of monogenic fields in $D^-$ which vanish
at infinity.
The Hardy subspaces $E^\pm L_2=\ran(E^\pm)$ are also the ranges of the Hardy projection
operators $E^\pm$ in $L_2(\Sigma)$ from Theorem~\ref{thm:cmcm}, which explains the notation.
There is a one-to-one correspondence between $f= F|_{\Sigma}\in E^\pm L_2$ and 
$F= C^\pm f: D^\pm\rightarrow \wedge\R^n$, and we sometimes identify $F$ and $f$, referring to
$F$ as belonging to the Hardy type subspace.

In the second splitting $L_2= N^+L_2\oplus N^-L_2$, which is pointwise, the subspace $N^+L_2$
consists of all fields tangential to $\Sigma$, and $N^-L_2$ consists of all fields normal
to $\Sigma$.

\begin{defn}   \label{defn:normtang}
  A multivector field $f: \Sigma\rightarrow \wedge\R^n$ is {\em tangential} if $\nu(x)\lctr f(x)=0$
  for almost all $x\in \Sigma$, and it is {\em normal} if $\nu(x)\wedg f(x)=0$ for almost all $x\in \Sigma$.
\end{defn}

The two projection operators $N^\pm$ are
$$
  N^+g:= \nu\lctr(\nu\wedg g)\qquad\text{and}\qquad N^-g:=\nu\wedg(\nu\lctr g).
$$
This tangential/normal splitting is orthogonal, so that $\|g_++g_-\|^2= \|g_+\|^2+\|g_-\|^2$,
when $g_\pm\in N^\pm L_2$.
On the other hand, the Hardy space splitting is not orthogonal, only topological
in the sense that $\|f_++f_-\| \approx \|f_+\|+\|f_-\|$ when $f_\pm\in E^\pm L_2$.

The operator theoretic problem underlying boundary value problems is to understand 
the relation between the splitting 
$E^+L_2\oplus E^-L_2$ and the splitting $N^+L_2\oplus N^-L_2$.
\begin{ex}   \label{ex:diracbvp}
  Consider the following basic Dirac BVP consisting in finding
$F: D^+\rightarrow \wedge\R^n$ solving the Dirac equation $(d+\del)F=0$ in $D^+$ with
given normal part $g$ on the boundary $\Sigma$. 
Under appropriate regularity assumptions, this means exactly that we are looking for
$f=F|_\Sigma \in E^+L_2$ such that $N^-f=g$.
Uniqueness and existence of such $f$, for each $g$, is clearly equivalent to
the restricted projection 
$$
  N^- :E^+L_2\longrightarrow N^-L_2
$$
being an isomorphism.
\end{ex}

In general, there are topological obstructions preventing $N^-|_{E^+L_2}$ from being an isomorphism.
However, modulo finite dimensional spaces, the operator is always invertible.

\begin{thm}    \label{thm:rellich}
  Let $\Sigma$ be any strongly Lipschitz surface. 
  Then the restricted projection $N^- :E^+L_2\rightarrow N^-L_2$ is a Fredholm
  operator of index $0$, i.e. has closed range and finite dimensional kernel and cokernel
  of equal dimensions.
  The same is true for all eight restricted projections 
\begin{alignat*}{2}  
  N^+ &: E^\pm L_2\longrightarrow N^+L_2, & \qquad   N^- &: E^\pm L_2\longrightarrow N^-L_2, \\
  E^+ &: N^\pm L_2\longrightarrow E^+L_2, & \qquad  E^- &: N^\pm L_2\longrightarrow E^-L_2.
\end{alignat*}
  If $\Sigma$ is a Lipschitz graph, then all these maps are isomorphisms.
\end{thm}
  
The key ingredient in the proof is a Rellich type estimate. 
The strong Lipschitz condition 
on $\Sigma$, i.e. that $\Sigma$ is locally the graph of a Lipschitz function, shows the 
existence of a smooth vector field $\theta$ which is transversal to $\Sigma$, i.e. $\scl{\nu(x)}{\theta(x)}\ge c>0$
for all $x\in\Sigma$.
Basic identities for the Clifford and interior products show that
$$
  |f|^2\scl\nu\theta = |\inv f|^2\tfrac 12(\nu\theta+\theta\nu)= \scl{\inv f\nu}{\inv f \theta}
  =\scl{-2\nu\lctr f+\nu f}{\inv f\theta},
$$
where $f\mapsto \inv f$ is the automorphism which negates $k$-vector fields with odd $k$.
Thus an application of Stokes' theorem yields the following Rellich identity
\begin{equation}    \label{eq:rellichid}
  \int_\Sigma |f|^2\scl\nu\theta \, d\sigma(y)=-2\int_\Sigma\scl{\nu\lctr f}{\inv f\theta}\, d\sigma(y)+ 
  \sum_{j=1}^n\iint_{D^+}\scl{F}{e_j\inv F(\pd_j\theta)} \, dx,
\end{equation}
for all $f=F|_\Sigma\in E^+L_2$, and therefore the estimate
$\|f\|\lesssim \| N^-f \|+\|F\|_{L_2(\supp\theta)}$.
If $\Sigma$ is a Lipschitz graph, we can choose $\theta=-e_n$, in which case the last term vanishes
and it follows that the restricted projection $N^- :E^+L_2\rightarrow N^-L_2$ is injective  
and has closed range. 
More generally, the map $f\mapsto F$ in the last term in the estimate can be shown to be 
compact, from which it follows that $N^- :E^+L_2\rightarrow N^-L_2$ has finite dimensional
null space and closed range. 
Finally the index of the restricted projection can be shown to be zero through either a duality argument
or the method of continuity. 
For details we refer to \cite{Ax2}.

As we shall make frequent reference to it, let us state the well known method of continuity.
\begin{thm}[Method of continuity]
  Let $\mX$ and $\mY$ be Banach spaces, and assume that $T_\lambda:\mX\rightarrow \mY$, $\lambda\in[0,1]$,
  is a family of bounded operators depending continuously on $\lambda$.
  If $T_\lambda$ are all semi-Fredholm operators, i.e. has closed range and finite dimensional null space,
  then the index, i.e. $\dim(\mY/\ran(T_\lambda))-\dim \nul(T_\lambda)$, of all operators $T_\lambda$ are equal.
\end{thm}
  
Another operator in $L_2(\Sigma)$ of importance to us, besides 
$E^\pm$ and $N^\pm$, is the following unbounded first order differential operator $\Gamma$.
\begin{defn}    \label{defn:gamma}
  Let $\Sigma$ be a strongly Lipschitz surface. 
  Denote by $\Gamma$ the unique closed operator in $L_2(\Sigma)$ with dense domain
  $\dom(\Gamma)\subset L_2(\Sigma)$ and the following action.
\begin{itemize}
\item[{\rm (i)}]
  If $f= f_++f_-= F_+|_\Sigma+ F_-|_\Sigma\in E^+L_2\oplus E^-L_2$ in the Hardy space
  splitting, then $\Gamma$ acts by exterior ($=-$ interior) differentiation on the monogenic fields $F^\pm$ 
in $D^\pm$ as
$$
  \Gamma f= (dF_+)|_\Sigma + (dF_-)|_\Sigma= (-\del F_+)|_\Sigma + (-\del F_-)|_\Sigma.
$$
\item[{\rm (ii)}]
  If $f= f_1+ \nu\wedg f_2\in N^+L_2\oplus N^-L_2$ in the tangential/normal
  splitting, where $f_1,f_2\in N^+L_2$, then $\Gamma$ acts by tangential exterior and interior differentiation 
  on the two parts respectively as
$$
  \Gamma f= d_\Sigma f_1 + \nu\wedg(\del_\Sigma f_2),
$$
where $d_\Sigma$ and $\del_\Sigma$ denote the intrinsic tangential exterior and interior differentiation operators
on the surface $\Sigma$.
\end{itemize}
\end{defn}
Recall that a bilipschitz parametrization $\rho:\R^{n-1}\rightarrow \Sigma$, locally around a point $y\in \Sigma$,
induces a pullback $\rho^*$, mapping tangential multivector fields $N^+L_2$ to $L_2(\R^{n-1};\wedge\R^{n-1})$.
Exterior differentiation commutes with this pullback, i.e. $d_\Sigma f= (\rho^*)^{-1}d_{\R^{n-1}}\rho^* f$.
Dual to this, a reduced pushforward $\tilde \rho_*$ in the terminology of \cite{Ax}, intertwines 
$\del_\Sigma$ and $\del_{\R^{n-1}}$.   

For Definition~\ref{defn:gamma} to make sense, one needs to show that the operators in (i) and (ii) coincide.
This is a consequence of the following proposition.
\begin{prop}   \label{prop:intertwining}
  The following intertwining relations for exterior and interior differentiation operators hold.
\begin{itemize}
\item[{\rm (i)}]
  If $U: D^\pm\rightarrow \wedge\R^n$ and $N_*(U), N_*(dU)\in L_2(\Sigma)$, then
 $N^+ U|_\Sigma\in\dom(d_\Sigma)$ and
$\Gamma (N^+ U|_\Sigma)=d_\Sigma(N^+ U|_\Sigma)=N^+ (dU|_\Sigma)$.
  If $U: D^\pm\rightarrow \wedge\R^n$ and $N_*(U)$, $N_*(\del U)\in L_2(\Sigma)$, then
  $\nu\lctr U|_\Sigma\in\dom(\del_\Sigma)$ and
$\Gamma (N^- U|_\Sigma)=\nu\wedg \del_\Sigma(\nu\lctr U|_\Sigma)=N^- (\del U|_\Sigma)$.
\item[{\rm (ii)}] 
  If $h\in N^+ L_2$, $d_\Sigma h\in L_2(\Sigma)$ and $x\notin\Sigma$, then
$$
  d\int_\Sigma E(y-x)\nu(y) h(y) \, d\sigma(y)= \int_\Sigma E(y-x)\nu(y) (d_\Sigma h)(y) \, d\sigma(y).
$$
  If $h\in N^- L_2$, $\del_\Sigma (\nu\lctr h)\in L_2(\Sigma)$ and $x\notin\Sigma$, then
$$
  d\int_\Sigma E(y-x)\nu(y) h(y) \, d\sigma(y)= \int_\Sigma E(y-x)\nu(y) (\nu\wedg \del_\Sigma (\nu\lctr h))(y)\, d\sigma(y).
$$
\end{itemize}
\end{prop}
The trace result $N^+ (dU|_\Sigma)= d_\Sigma(N^+ U|_\Sigma)$ in (i) is a special case
of the fundamental fact that the exterior differentiation and pullbacks commute.
Indeed, if $i: \Sigma\rightarrow \R^n$ denotes inclusion, then $N^+ U|_\Sigma= i^*(U)$.
The trace result $N^- (\del U|_\Sigma)= \nu\wedg \del_\Sigma(\nu\lctr U|_\Sigma)$
can then be obtained by Hodge star duality.
For more details of the proof of Proposition~\ref{prop:intertwining}, we refer to
\cite[Proposition 4.10]{Ax2} and \cite[Proposition 6.2.5]{Ax}.

The usefulness and relevance of the operator $\Gamma$ to this paper, is that it is closely 
related to the decomposition (\ref{eq:homoparts}), as the following proposition shows.
\begin{prop}   \label{prop:gammarelevance}
The operator $\Gamma$ relates to the splitting of $\wedge\R^n$ into homogeneous $k$-vectors as follows.
\begin{itemize}
\item[{\rm (i)}]
  If $h\in L_2(\Sigma;\wedge^k)$ and $\Gamma(N^+ h)=0$, then $(C^\pm h)_{k+2}=0$.
  If $h\in L_2(\Sigma;\wedge^k)$ and $\Gamma(N^- h)=0$, then $(C^\pm h)_{k-2}=0$.
  Thus, if $h\in L_2(\Sigma;\wedge^k)$ and $\Gamma h=0$, then $C^\pm h$ is a $k$-vector field.
\item[{\rm (ii)}]
 For a monogenic field $F=\sum_{k=0}^n F_k$ in $D^\pm$, i.e. $(d+\del)F=0$, 
 where $F_k:D^\pm\rightarrow \wedge^k\R^n$, 
 the following are equivalent.
 \begin{itemize}
\item All homogeneous parts $F_k$ of $F$ are monogenic, i.e. $(d+\del)F_k=0$.
\item $F$ is two-sided monogenic, i.e. $\dirac F=0= F\dirac$, where
$F\dirac:= \sum_{j=1}^n (\pd_j F)e_j$.
\item $F$ satisfies $dF=\del F=0$.
\end{itemize}
\end{itemize}
\end{prop}
The result (i) shows that $\Gamma$ governs the off-diagonal mapping of the Cauchy integral,
i.e. the first and last terms in (\ref{eq:Cauchymapping}). 
Indeed, an integration by parts rewrites these terms as single layer potentials
\begin{align*}
  \int_\Sigma E(y-x)\wedg\nu(y)\wedg h(y) \, d\sigma(y) & = \int_\Sigma \Phi(y-x)\, \nu(y)\wedg (\Gamma h(y)) \, d\sigma(y), \qquad x\notin\Sigma, \\
  \int_\Sigma E(y-x)\lctr(\nu(y)\lctr h(y))\, d\sigma(y) & = -\int_\Sigma \Phi(y-x)\, \nu(y)\lctr (\Gamma h(y)) \, d\sigma(y), \qquad x\notin\Sigma,
\end{align*}
where $\Phi(x)$ denotes the fundamental solution for the Laplace operator and $E(x)=\nabla\Phi(x)$.
For proof of Proposition~\ref{prop:gammarelevance}, 
we refer to \cite[Lemma 4.13, Proposition 4.5]{Ax2}.

Since $\Gamma$ acts by exterior and interior differentiation, it is clear that $\Gamma^2=0$, or more precisely
$\ran(\Gamma)\subset \nul(\Gamma)$.
We have inclusions of function spaces
$$
  L_2^R(\Sigma)\subset L_2^N(\Sigma)\subset L_2^D(\Sigma)\subset L_2(\Sigma),
$$
where $L_2^R(\Sigma):= \ran(\Gamma)$, $L_2^N(\Sigma):= \nul(\Gamma)$
and $L_2^D(\Sigma):=\dom(\Gamma)$.
Here $L_2^N(\Sigma)$ is always a closed subspace of $L_2(\Sigma)$ and $L_2^D(\Sigma)$ is always a Hilbert 
space densely embedded in $L_2(\Sigma)$.
The domain $L_2^D(\Sigma)$ is equipped with the graph norm $\|f\|_D^2:= \|f\|^2_2+ \|\Gamma f\|^2_2$,
which makes it a Hilbert space.
The range $L_2^R(\Sigma)$ is equipped with the range norm
$$
  \|f\|_R^2:= \inf\sett{\|\Gamma u\|_2^2+ \|u\|_2^2}{u\in\dom(\Gamma), \Gamma u=f},
$$
which makes it a Hilbert space.
The properties of the range $L_2^R(\Sigma)$ depends on the surface $\Sigma$, as the following lemma
shows.
\begin{lem}    \label{lem:cohom}
  If $\Sigma$ is an unbounded Lipschitz graph, then $L_2^R(\Sigma)$ is dense and not closed in $L_2^N(\Sigma)$.
  If $\Sigma$ is a bounded Lipschitz surface, then $L_2^R(\Sigma)$ is a closed subspace of $L_2^N(\Sigma)$
  of finite codimension.
  In particular, if $D^+$ is Lipschitz diffeomorphic to the unit ball, then the codimension is $4$ and
$$
  L_2^R(\Sigma)=\sett{f\in L_2^N(\Sigma)}{f_0=f_n= \int_\Sigma \nu\wedg f_{n-1} d\sigma= \int_\Sigma \nu\lctr f_1 d\sigma=0},
$$
where $f_k:\Sigma\rightarrow \wedge^k \R^n$ denotes the $k$-vector part of $f$.
\end{lem}

\begin{proof}
Since  $N^+L_2$ and $N^-L_2$ are invariant under $\Gamma$,
we may consider tangential and normal multivector fields separately.
Moreover, the two operators $d$ and $\del$ acting in $\R^n$ satisfy
$$
  \del(F e_\nindex) = (d F) e_\nindex,
$$ 
and taking the normal part of the trace of this identity, 
Proposition~\ref{prop:intertwining}(i) shows that
$\Gamma(fe_\nindex)= (\Gamma f)e_\nindex$ for all tangential $f\in N^+ L_2$.
Thus the actions of $\Gamma$ on $N^+L_2$ and $N^-L_2$ are similar,
and it suffices to consider $\Gamma=d_\Sigma$ acting on tangential multivector field.
Here the stated results are well known facts from de Rham cohomology.
\end{proof}

Instead of working with the projections $E^\pm$ and $N^\pm$, it is often convenient to work
with the associated reflection operators $E:= E^+-E^-$ and $N:= N^+-N^-$, where $E^2=N^2=I$,
which explains the naming. 
Here $E$ is the principal 
value Cauchy singular integral from Theorem~\ref{thm:cmcm} and 
$$
  Nf= \nu\inv f \nu
$$
is the operator which reflects a multivector field across $\Sigma$.
Following the boundary equation method developed in \cite{Ax, Ax1}, we shall make use of the
{\em rotation operator}
$$
  ENf(x)= \pv \int_\Sigma E(y-x) \inv f(y) \nu(y) d\sigma(y).
$$
The important connection between $EN$ and the restricted projections above, is that
\begin{alignat*}{2}
  I+EN &= 2(E^+N^+ + E^-N^-), &\qquad   N(I+EN)N &= 2(N^+E^+ + N^-E^-),  \\
  I-EN &= 2(E^+N^-+ E^-N^+),  &\qquad  N(I-EN)N &= 2(N^+E^-+ N^-E^+).
\end{alignat*}
For example, this shows that $I+EN$ is the direct sum of the restricted projections
$E^+:N^+L_2\rightarrow E^+L_2$ and $E^-: N^-L_2\rightarrow E^-L_2$.
Thus, in order to prove that all eight restricted projections are Fredholm operators,
it suffices to prove that the two operators $I\pm EN$ are Fredholm operators on the
full space $L_2(\Sigma;\wedge)$.
We record the following generalization of Theorem~\ref{thm:rellich}, which was proved
in \cite[Theorem 4.15]{Ax2}
through Rellich estimates involving a pair of monogenic fields $F^\pm :D^\pm \rightarrow \wedge\R^n$
and the method of continuity.
\begin{thm}    \label{thm:rotresolvent}
  Let $\Sigma$ be a strongly Lipschitz surface. Then $\lambda+ EN: L_2(\Sigma)\rightarrow L_2(\Sigma)$ 
  is a Fredholm operator
  with index zero for all $\lambda\in \R$ (and more generally in a double sector around the real axis). 
\end{thm}  

The last result we shall need is the following analogue of Theorem~\ref{thm:rellich} for the
subspaces $L_2^R(\Sigma)$, $L_2^N(\Sigma)$ and $L_2^D(\Sigma)$.
\begin{thm}    \label{thm:domainrellich}
  All four projections $E^\pm$ and $N^\pm$ leave each of the subspaces $L_2^R(\Sigma)$,
  $L_2^N(\Sigma)$ and $L_2^D(\Sigma)$ invariant and act boundedly in them.
  All eight restricted projections 
\begin{alignat*}{2}  
  N^+ &: E^\pm L_2^x\longrightarrow N^+L_2^x, & \qquad   N^- &: E^\pm L_2^x\longrightarrow N^-L_2^x, \\
  E^+ &: N^\pm L_2^x\longrightarrow E^+L_2^x, & \qquad  E^- &: N^\pm L_2^x\longrightarrow E^-L_2^x
\end{alignat*}
are Fredholm operators, for $x=R,N,D$.
All eight maps are injective when $x=R$, i.e. when acting in the range $L_2^R(\Sigma)$,
for all strongly Lipschitz surfaces $\Sigma$.
\end{thm}

\begin{proof}[Proof of Theorem~\ref{thm:domainrellich}]
  That $E^\pm$ and $N^\pm$ act boundedly in all three subspaces is clear from (i) and (ii) in 
  Definition~\ref{defn:gamma}.
  As noted above, the Fredholm property of all eight restricted projections will follow if we
  prove that $I\pm EN$ are Fredholm operators on $L_2^x(\Sigma)$.

(1)
  Fredholmness of the operators acting in $L_2^N(\Sigma)$ and 
  $L_2^D(\Sigma)$ follows from Theorem~\ref{thm:rotresolvent} and the method of continuity.
  For details, we refer to \cite[Theorem 4.15]{Ax2} where it was shown that 
  $I\pm EN$ are Fredholm operators
  with index zero on $L_2(\Sigma)$, $L_2^D(\Sigma)$ and $L_2^N(\Sigma)$.

  To show that $I\pm EN$ are Fredholm operators on $L_2^R(\Sigma)$, note that this is 
  immediate from Lemma~\ref{lem:cohom} and the result for $L_2^N(\Sigma)$ when 
  $\Sigma$ is bounded. In case of a 
  Lipschitz graph, consider the commutative diagram
$$
\xymatrix@1{ 
  0 \ar[r]  & L_2^N(\Sigma) \ar[r]^i \ar[d]^{I\pm EN} & L_2^D(\Sigma) \ar[r]^\Gamma \ar[d]^{I\pm EN} & L_2^R(\Sigma) \ar[r]  \ar[d]^{I\pm EN} & 0 \\
  0 \ar[r]  & L_2^N(\Sigma) \ar[r]^i & L_2^D(\Sigma) \ar[r]^\Gamma & L_2^R(\Sigma) \ar[r] & 0. 
}
$$
Note that the rows are exact, i.e. the inclusion map $i$ is injective, $\Gamma$ is surjective 
and $\nul(\Gamma)=\ran(i)= L_2^N(\Sigma)$.
It has been shown that the first two vertical maps are Fredholm operators. 
We can now apply a general technique, the five lemma, to deduce that 
$I\pm EN:L_2^R(\Sigma)\rightarrow L_2^R(\Sigma)$ is a Fredholm operator as well.
For details concerning the five lemma for Fredholm operators we refer to Pryde~\cite{P}.

(2)
It remains to show injectivity on $L_2^R(\Sigma)$. For this, it suffices to show that 
$E^\pm L_2^R \cap N^\pm L_2^R=\{0\}$  for all four intersections.
If $\Sigma$ is a Lipschitz graph, then the result follows from Theorem~\ref{thm:rellich}
since in this case $E^\pm L_2\cap N^\pm L_2=\{0\}$.
Assuming that $\Sigma$ is a bounded Lipschitz surface,
consider for example $E^-L_2^R\cap N^+L_2^R$, and let $F=dU=-\del U$, where 
$F$ and $U$ are monogenic in $D^-$, $N_*(U), N_*(F)\in L_2(\Sigma)$, $U,F\rightarrow 0$
when $x\rightarrow \infty$, and $\nu\lctr F=0$ on $\Sigma$.
An application of Stokes' theorem shows that
$$
  \iint_{D^-}|F|^2 \, dx= \int_\Sigma\scl{U}{\nu\lctr F} \, d\sigma(y)-\iint_{D^-}\scl{U}{\del F} \, dx=0,
$$
and therefore $F=0$.
Note that the asymptotics of the Cauchy kernel $E(x)$ shows that $|U|\lesssim |x|^{1-n}$,
$|F|\lesssim |x|^{-n}$ and $|\del F|\lesssim |x|^{-1-n}$ as $x\rightarrow\infty$,
so that both $D^-$ integrals are convergent.

A similar argument shows that the other three intersections also equal $\{0\}$.
\end{proof}

\begin{cor}   \label{cor:infdim}
  If $1\le k\le n-1$, then the Hardy space $E^\pm L_2^N(\Sigma;\wedge^k)$ of boundary traces of 
  two-sided monogenic $k$-vector fields in $D^\pm$ is infinite dimensional.
\end{cor}
\begin{proof}
  Consider the Fredholm operator
$$
  E^\pm : N^+ L_2^N\longrightarrow E^\pm L_2^N.
$$
By Proposition~\ref{prop:gammarelevance}, this maps $k$-vector fields to $k$-vector fields.
Since the space of tangential $f\in N^+L_2(\Sigma;\wedge^k)$ such that $d_\Sigma f=0$ 
is infinite dimensional, the corollary follows.
\end{proof}

\section{Hilbert transforms for $k$-vector fields}    \label{sec:kHilbert}

In this section we prove the following main result of this paper, which extends
Theorem~\ref{thm:main} to more general $k$-vector fields.
\begin{thm}   \label{thm:maink}
  Let $D\subset \R^n$ be a Lipschitz graph, interior or exterior domain and assume that 
  $1\le k\le n-1$.

  Let $U:D \rightarrow\wedge^k\R^n$ be such that 
  $\del dU=0= d\del U$ and 
  $N_*(U), N_*(dU), N_*(\del U)\in L_2(\Sigma)$.
  If $D$ is an exterior domain, also assume that 
  $U, dU, \del U\rightarrow 0$ when $x\rightarrow \infty$.
  If $D$ is an interior domain, also assume that $\del U=0$ if $k=1$ and $dU=0$ if $k=n-1$.
  Then there are unique Cauchy type harmonic conjugates $V_1:D\rightarrow \wedge^{k-2}\R^n$
  and $V_2: D\rightarrow \wedge^{k+2}\R^n$ such that
$$
  \|N_*(V_1)\|_2+ \|N_*(V_2)\|_2 \lesssim \|N_*(U)\|_2+ \|N_*(dU)\|_2+ \|N_*(\del U)\|_2. 
$$
\end{thm}

Before the proof of the theorem, we make some remarks.
In the scalar case, we could apply the known results from Theorem~\ref{thm:dahlbergkenig} 
on $L_p$-invertibility of the classical double layer potential operator to  
prove existence and uniqueness of Cauchy type conjugate functions in Theorem~\ref{thm:main}.
For $k$-vector fields, $1\le k\le n-1$, Theorem~\ref{thm:dahlbergkenig} is no longer
available. In particular, we do not obtain any $L_p$-results for $p>2$, since the atomic
estimates in the proof of Theorem~\ref{thm:dahlbergkenig} in an essential way use
that the equation is scalar.
Instead, we make use of a natural $L_2$-based boundary function space, $L_2^D(\Sigma)$, of mixed $0$ and $1$
order regularity. The key observation is that both $dU$ and $\del U$ are (two-sided) monogenic when
$U$ is a two-sided harmonic $k$-vector field. Thus, in order to apply the well established
$L_2$-theory for BVPs, we need to require that $N_*(U)$, $N_*(dU)$ and $N_*(\del U)$ belong to $L_2$.

Just like in the scalar case, the Cauchy type harmonic conjugate functions to a $k$-vector field
$U$ can be calculated using a generalized double layer potential operator.
Indeed, according to (\ref{eq:Cauchymapping}), the Cauchy integral maps
$$
\xymatrix{ 
  & & V_2 \\
  h \ar[rru]^>>>>{(C^\pm)_{k+2}} \ar[rr]^>>>>{(C^\pm)_k} \ar[rrd]_>>>>{(C^\pm)_{k-2}} & & U \\
  & & V_1.
 }
$$
Thus, we need to solve for $h$ in the generalized double layer potential equation
$(C^\pm h)_k=U$.

\begin{cor}   \label{cor:gendlp}
  Let $D^\pm\subset\R^n$ be a Lipschitz graph, interior or exterior domain, and let $1\le k\le n-1$.
  Then the range and null space of $(C^\pm)_k$, with domain $L_2^D(\Sigma;\wedge^k)$, are
\begin{align*}
  \ran((C^\pm)_k) &= \sett{U: D^\pm\rightarrow \wedge^k\R^n}{\del dU=0=d\del U, N_*(U), N_*(\del U), N_*(dU)\in L_2(\Sigma)}, \\
  \nul((C^\pm)_k) &= \sett{F|_{\Sigma}}{F: D^\mp\rightarrow \wedge^k\R^n,
     dF=0=\del F, N_*(F)\in L_2(\Sigma)},
\end{align*}
  with the same modifications of $\ran((C^\pm)_k)$ as in Theorem~\ref{thm:maink} when $D^\pm$ is 
  an exterior domain and when $D^\pm$ is an interior domain and $k=1,n-1$,
  and where $F\to 0$ when $x\to\infty$ when $D^\mp$ is an exterior domain and $F\in\nul((C^\pm)_k)$.
  The operator
$$
  (C^\pm)_k: L_2^D(\Sigma;\wedge^k) /\nul((C^\pm)_k)\longrightarrow \ran((C^\pm)_k)
$$
is an isomorphism.
Thus, if $U\in \ran((C^\pm)_k)$, its Cauchy type harmonic conjugates are
$$
  V_1= (C^\pm ((C^\pm)_k)^{-1} U)_{k-2}\qquad \text{and}\qquad 
  V_2= (C^\pm ((C^\pm)_k)^{-1} U)_{k+2}.
$$
\end{cor}

\begin{proof}
  If $U\in \ran((C^\pm)_k)$, i.e. if $U=(C^\pm h)_k$ for some dipole density $h\in L_2^D(\Sigma;\wedge^k)$,
  then it follows that $U$ is harmonic and $N_*(U), N_*(\del U), N_*(dU)\in L_2(\Sigma)$.
  But $\del dU= -\del \del (C^\pm h)_{k+2}=0$, so $U$ is in fact two-sided harmonic. 
  The converse inclusion follows from the existence proof of Theorem~\ref{thm:maink} below.

  Clearly, if $F: D^\mp\rightarrow \wedge^k\R^n$ satisfies $dF=0=\del F$ and $N_*(F)\in L_2(\Sigma)$,
  then $F|_\Sigma$ belong to the Hardy subspace for the complementary domain $D^\mp$.
  In particular $C^\pm (F|_\Sigma)=0$ and $F|_\Sigma\in \nul((C^\pm)_k)$.
  The converse inclusion follows from the uniqueness proof of Theorem~\ref{thm:maink} below.
\end{proof}

\begin{rem}   \label{rem:compressions}
  According to (\ref{eq:dlpdefn}), the classical double layer potential is the compression of 
  the Cauchy integral to the subspace of scalar valued functions.
  In Corollary~\ref{cor:gendlp}, identifying a harmonic function with its boundary trace,
  we have generalized this by compressing the Cauchy integral/Hardy projection $E^\pm$
  to the operator $(E^\pm)_k$, acting in the subspace of $k$-vector fields.
  Other useful compressions of the Cauchy integral use instead the subspaces of tangential
  or normal multivector fields $N^\pm L_2(\Sigma)$, which are relevant for BVPs.
  For example, consider the BVP in Example~\ref{ex:diracbvp} consisting in finding $F:D^+\to \wedge \R^n$
  satisfying the Hodge--Dirac equation, with a prescribed normal part $g$ of the trace $f=F|_\Sigma$.
  Equivalently, we look for a Hardy function $f\in E^+L_2$ satisfying $N^-f=g$.
  Making the ansatz $f= E^+ h$, with $h\in N^- L_2$, we obtain instead a double layer type
  equation $N^-E^+ h=g$ in the subspace $N^-L_2$.
  As shown in \cite{Ax1}, the well posedness of the BVP is essentially equivalent to
  the compressed Cauchy integral $N^-E^+|_{N^-L_2}$ being an isomorphism.
  We note that these types of compressions to $N^\pm L_2$ in general have better properties than
  the compressions $(E^\pm)_k|_{L_2(\Sigma;\wedge^k)}$ used in Corollary~\ref{cor:gendlp},
  as $N^\pm E^\pm|_{N^\pm L_2}$ are Fredholm operators.
\end{rem}

\begin{proof}[Uniqueness proof of Theorem~\ref{thm:maink}]
Assume that $h\in L_2(\Sigma;\wedge^k)$ is such that its Cauchy extension satisfies
$U=(C^\pm h)_k=0$. We aim to prove that $h\in L_2^N(\Sigma)$ and $C^\pm h=0$,
so that the Cauchy type harmonic conjugates $V_1=(C^\pm h)_{k-2}$ and
$V_2=(C^\pm h)_{k+2}$ vanish.

Define monogenic fields
\begin{align*}
  \tilde V_1(x) & := C^\pm (N^- h)(x)= \pm\int_\Sigma E(y-x)(\nu(y)\lctr h(y)) \, d\sigma(y), \\
  \tilde V_2(x) & := C^\pm (N^+ h)(x)= \pm\int_\Sigma E(y-x)(\nu(y)\wedg h(y))\, d\sigma(y),
\end{align*}
for $x\in D^\pm$,
so that $V_1= (\tilde V_1)_{k-2}$, $V_2=(\tilde V_2)_{k+2}$ and
$(\tilde V_1)_k+ (\tilde V_2)_k=U=0$ by assumption.
It follows that $dV_2= ((d+\del)\tilde V_2)_{k+3}=0$ and
$\del V_2= -d(\tilde V_2)_k= d(\tilde V_1)_{k}=((d+\del)\tilde V_1)_{k+1}=0$.
Thus $V_2$ and therefore $\tilde V_2$ are two-sided monogenic and
\begin{equation}   \label{eq:gammaiszero}
  \Gamma (E^\pm N^+h)=0.
\end{equation}
Similarly, it follows that $\Gamma(E^\pm N^-h)=0$.

We first show that (\ref{eq:gammaiszero}) implies that $h_1:=N^+h$ has regularity 
$h_1\in L_2^D(\Sigma)$.
Note that $E^\pm h_1= \tfrac 12(I\pm E)h_1= \tfrac 12(I\pm EN)h_1$, and
consider the commutative diagram
$$
\xymatrix@1{ 
  L_2^D(\Sigma) \ar[r]^{I\pm EN} \ar[d]^i & L_2^D(\Sigma) \ar[d]^i \\
  L_2(\Sigma) \ar[r]^{I\pm EN} & L_2(\Sigma), 
}
$$
where the inclusion $i$ is dense. 
Moreover, as explained in the proof of Theorem~\ref{thm:domainrellich}, 
the method of continuity shows that $I\pm EN$
are Fredholm operators with index zero on both $L_2(\Sigma)$ and 
$L_2^D(\Sigma)$.
Since $(I\pm EN) h_1\in L_2^D(\Sigma)$, a general regularity theorem
for Fredholm operators \cite[Proposition 4.16]{Ax2} shows that $h_1\in L_2^D(\Sigma)$.
Similarly, $N^- h\in L_2^D(\Sigma)$.

We have shown that $E^\pm (\Gamma N^+h)=0= E^\pm (\Gamma N^- h)$,
i.e. $\Gamma N^+ h\in E^\mp L_2^R\cap N^+ L_2^R$ and
$\Gamma N^- h\in E^\mp L_2^R\cap N^- L_2^R$.
Thus Theorem~\ref{thm:domainrellich} shows that $\Gamma(N^+ h)= \Gamma(N^-h)=0$,
so that $h\in L_2^N(\Sigma)$.
According to Proposition~\ref{prop:gammarelevance}(i), we have 
$V_1=V_2=0$ and thus $C^\pm h=0$.
\end{proof}

\begin{proof}[Existence proof of Theorem~\ref{thm:maink}]
  Assume that $U:D\rightarrow \wedge^k\R^n$ is twosided harmonic in the sense that
  $\del dU=0= d\del U$, where $N_*(U), N_*(dU), N_*(\del U)\in L_2(\Sigma)$.
  If $D$ is an exterior domain, also assume that 
  $U, dU, \del U\rightarrow 0$ when $x\rightarrow \infty$.
  If $D$ is an interior domain, also assume that $\del U=0$ if $k=1$ and $dU=0$ if $k=n-1$.

  We aim to construct $V_1: D\rightarrow \wedge^{k-2}\R^n$ and $V_2:D\rightarrow \wedge^{k+2}\R^n$
  and a dipole density $h\in L_2^D(\Sigma;\wedge^k)$ such that $C^\pm h= V_1+U+V_2$
  and 
$$
  \|N_*(V_1)\|_2+ \|N_*(V_2)\|_2+\|h\|_D \lesssim \|N_*(U)\|_2+ \|N_*(dU)\|_2+ \|N_*(\del U)\|_2. 
$$

(1)
We first construct a tangential $k$-vector field $h_2\in N^+ L_2^D(\Sigma;\wedge^k)$ such that
$dU= d(C^\pm h_2)$.
To this end, consider the singular integral equation
$$
  N^+ E^\pm \tilde h_2= N^+ (dU|_\Sigma). 
$$
From the assumption, $dU$ is a monogenic field, and Proposition~\ref{prop:intertwining}(i)
shows that $N^+(dU|_\Sigma)\in L_2^R(\Sigma)$.
We claim that the compressed Cauchy integral
$$
  N^+ E^\pm: N^+L_2^R(\Sigma)\longrightarrow N^+ L_2^R(\Sigma)
$$
is an isomorphism. Since it is the composition of $E^\pm: N^+ L_2^R(\Sigma)\rightarrow E^\pm L_2^R(\Sigma)$
and $N^+: E^\pm L_2^R(\Sigma)\rightarrow N^+L_2^R(\Sigma)$, it follows from 
Theorem~\ref{thm:domainrellich} that it is an injective Fredholm operator.
Using the operator algebra developed for boundary value problems in \cite{Ax1}, we see that
\begin{align*}
  \lambda^2 -4(N^+E^+N^+ + N^-E^-N^-) &= (\lambda N)^2- (N+E)^2 \\
  &= (\lambda+1+EN)N(\lambda-1-EN)N, \\
  \lambda^2 -4(N^+E^-N^+ + N^-E^+N^-) &= (\lambda N)^2- (N-E)^2 \\
  &= (\lambda+1-EN)N(\lambda-1+EN)N.
\end{align*}
The right hand sides are seen to be Fredholm operators in $L_2^R(\Sigma)$ for all real $\lambda$ as in 
part (1) in the proof of Theorem~\ref{thm:domainrellich}.
Applying the left hand sides to $f\in N^+ L_2^R(\Sigma)$, so that $N^-f=0$,
shows that $\lambda^2 -N^+ E^\pm: N^+L_2^R(\Sigma)\rightarrow N^+ L_2^R(\Sigma)$
are all Fredholm operators.
Since this operator clearly is invertible for large enough $\lambda$, the method
of continuity shows that $N^+ E^\pm: N^+L_2^R(\Sigma)\rightarrow N^+ L_2^R(\Sigma)$
has index zero, and therefore is surjective, since it has been shown to be injective.

Solving the equation, we obtain a unique $\tilde h_2\in N^+ L_2^R(\Sigma)$ such that
$$
  N^+(E^\pm \tilde h_2- (dU)|_\Sigma)=0,
$$
where we verify that $(dU)|_\Sigma\in E^\pm L_2^N(\Sigma)$.
It is here we need the topological assumption on $\Sigma$.
If $\Sigma$ is an unbounded graph, then $N^- L_2\cap E^\pm L_2=\{0\}$.
On the other hand, if $D^+$ is Lipschitz diffeomorphic to the unit ball and $1\le k\le n-3$, 
then $(dU)|_{\Sigma}\in E^\pm L_2^R(\Sigma)$ according to Lemma~\ref{lem:cohom}.
If $k=n-1$ then $dU=0\in E^\pm L_2^R(\Sigma)$ by assumption.
Finally, if $k=n-2$ then Stokes' theorem shows that
$$
  \int_\Sigma \nu\wedg (dU) d\sigma=\pm \int_{D^\pm} d(dU) dx=0
$$
since $d^2=0$. (Note that Stokes' theorem is applicable in $D^-$ since $|dU|\lesssim |x|^{-n}$
as $x\rightarrow \infty$.)
Thus $E^\pm \tilde h_2- (dU)|_\Sigma\in E^\pm L_2^R\cap N^- L_2^R$, which equals
$\{0\}$ by Theorem~\ref{thm:domainrellich}.
In either case, we conclude that 
$$
  dU= d(C^\pm h_2), \qquad\text{for some } h_2\in N^+L_2^D(\Sigma),
$$
such that $\tilde h_2= \Gamma h_2$ and 
$\|N_*(U)\|_2+\|N_*(dU)\|_2\gtrsim \|\tilde h_2\|_R\approx \|\tilde h_2\|_2+ \|h_2\|_2$.

A similar argument proves that
$$
  \del U= \del(C^\pm h_1), \qquad\text{for some } h_1\in N^-L_2^D(\Sigma),
$$
such that $\|h_1\|_2+ \|\Gamma h_1\|_2\lesssim \|N_*(U)\|_2+ \|N_*(\del U)\|_2$.

(2)
To construct Cauchy type harmonic conjugates $V_1$ and $V_2$, write
$\tilde V_i:= C^\pm h_i$, $i=1,2$.
Then
\begin{align*}
  d(U-(\tilde V_1)_k- (\tilde V_2)_k)= dU-0-dU=0, \\
  \del(U-(\tilde V_1)_k- (\tilde V_2)_k)= \del U-\del U-0=0,
\end{align*}
so that $(U-(\tilde V_1)_k- (\tilde V_2)_k)|_\Sigma\in E^\pm L^N_2$.
Thus, defining a dipole density
$$
  h:= (U-(\tilde V_1)_k- (\tilde V_2)_k)|_\Sigma + h_1+h_2\in L_2^D(\Sigma;\wedge^k),
$$
$V_1:= (\tilde V_1)_{k-2}$ and $V_2:=(\tilde V_2)_{k+2}$
gives
$$
  C^\pm h= U-(\tilde V_1)_k- (\tilde V_2)_k+ \tilde V_1+ \tilde V_2= V_1+U+V_2.
$$
Since $\|N_*(V_1)\|_2+ \|N_*(V_2)\|_2\lesssim \|N_*(U)\|_2+
\|h_1\|_2+ \|h_2\|_2\lesssim \|N_*(U)\|_2+ \|N_*(dU)\|_2+ \|N_*(\del U)\|_2$
and
$\|h\|_D\lesssim \|(U-(\tilde V_1)_k- (\tilde V_2)_k)|_\Sigma\|_2+ \|h_1\|_D+\|h_2\|_D 
\lesssim \|N_*(U)\|_2+ \|N_*(dU)\|_2+ \|N_*(\del U)\|_2$,
the proof of Theorem~\ref{thm:maink} is complete.
\end{proof}

\section{Other types of harmonic conjugates}    \label{sec:hodge}

Recall from the discussion in the introduction that if $U:D\rightarrow \wedge^k \R^n$ is a
two-sided harmonic $k$-vector field, two fields $V_1:D\rightarrow \wedge^{k-2}\R^n$
and $V_2:D\rightarrow \wedge^{k+2} \R^n$ are said to be conjugate to $U$ if
$(d+\del)(V_1+U+V_2)=0$.
As noted, further conditions need to be imposed on $V_1$, $V_2$ for this problem
to be well-posed, i.e. for $V_i$ to be (essentially) unique.
Theorems~\ref{thm:main} and \ref{thm:maink} show that the problem becomes well posed if $V_i$ are required
to be Cauchy type harmonic conjugates.
The following proposition expresses this condition as a boundary condition.

\begin{prop}
  Assume that $V_1$ and $V_2$ are harmonic conjugates to $U$ in $D=D^+$,
  with $N_*(U), N_*(V_1), N_*(V_2)\in L_2(\Sigma)$.
  Then they are of Cauchy type if and only if there exists $U^-:D^-\rightarrow\wedge^{k}\R^n$
  with harmonic conjugates $V_1^-:D^-\rightarrow\wedge^{k-2}\R^n$
  and $V_2^-:D^-\rightarrow\wedge^{k+2}\R^n$ in $D^-$, 
  with $N_*(U^-), N_*(V_1^-), N_*(V_2^-)\in L_2(\Sigma)$
  (and decay at infinity when $D^-$ is an exterior domain), such that
$$
  V_1|_\Sigma+ V_1^-|_\Sigma=0\qquad\text{and}\qquad V_2|_\Sigma+V_2^-|_\Sigma=0.
$$
\end{prop}
\begin{proof}
  Recall that by definition, $V_i$ are of Cauchy type if there exists $h:\Sigma\rightarrow \wedge^k\R^n$
  such that $C^+ h= V_1+U+V_2$.
  In this case, $V_1^-+U^-+V_2^-:= C^-h$ has the required properties since $E^+h+E^-h=h$.
  Conversely, if $V_1^-+U^-+V_2^-$ has the required properties, let
  $h:=U|_\Sigma + U^-|_\Sigma= (V_1+U+V_2)|_\Sigma +(V_1^-+U^-+V_2^-)|_\Sigma:\Sigma\rightarrow \wedge^k \R^n$.
  Then $C^+h= V_1+U+V_2$.
\end{proof}
Note that since Cauchy type conjugates are defined in terms of the Cauchy integral, which concerns
the interplay between monogenic fields in $D^+$ and $D^-$, the boundary condition above
is a {\em transmission problem}, i.e. a jump relation between pairs of monogenic fields in $D^\pm$.

We end this paper with a construction of harmonic conjugates which differ from the Cauchy type 
ones in general. To avoid technicalities, we shall only consider interior domains, i.e. 
$D$ is assumed to be Lipschitz diffeomorphic to the unit ball. 
(Unlike the situation for the Cauchy type conjugates, this second construction does not involve
the complementary domain $D^-$.)

Consider the exterior and interior derivative operators $d$ and $\del$ in (\ref{eq:complex}).
These are formally (anti-) adjoint in the sense that
$\int_{\R^n}\scl{\del F(x)}{G(x)}dx= -\int_{\R^n}\scl{F(x)}{d G(x)}dx$ for all $F, G\in C^\infty_0(\R^n;\wedge)$,
which is a consequence of the definition (\ref{eq:intprod}) of the interior product.
On the domain $D$, on the other hand, an application of Stokes' theorem shows that
$$
   \int_{D}\scl{\del F}{G}dx+ \int_{D}\scl{F}{d G}dx= \int_\Sigma\scl{\nu\lctr F}{G} d\sigma
   = \int_\Sigma\scl{F}{\nu\wedg G}d\sigma.
$$
Thus, in order to make $d$ and $-\del$ to be adjoint operators in $L_2(D;\wedge)$, one
needs to impose either tangential boundary conditions on $\del$ (so that $\nu\lctr F|_\Sigma=0$)
or normal boundary conditions on $d$ (so that $\nu\wedg G|_\Sigma=0$).
If these boundary conditions are imposed on the operators in a suitable weak sense,
then one obtains two pairs of densely defined, closed and adjoint operators $(-\udel, d)$ and $(-\del, \ud)$ 
in $L_2(D;\wedge)$, where the domains of the operators are
\begin{alignat*}{2}
  \dom(\del) &= \sett{F\in L_2}{\del F\in L_2}, &\qquad 
   \dom(\udel)&= \sett{F\in L_2}{\del F\in L_2, \nu\lctr F|_\Sigma=0}, \\
  \dom(d)&= \sett{G\in L_2}{d G\in L_2},&\qquad 
   \dom(\ud)&= \sett{G\in L_2}{d G\in L_2, \nu\wedg G|_\Sigma=0}, 
\end{alignat*}
i.e. $\udel$ is $\del$ with tangential boundary conditions and $\ud$ is $d$ with normal boundary conditions.
For further details, we refer to \cite[Section 4]{AMc}.
The following proposition summarizes relevant facts from the theory of Hodge decompositions 
that we need below to construct harmonic conjugates.
\begin{prop}   \label{prop:hodge}
  Assume that $D\subset \R^n$ is Lipschitz diffeomorphic to the unit ball.
\begin{itemize}
\item[{\rm (i)}]
  If $F\in L_2(D;\wedge^k)$ satisfies $d F=0$, and $2\le k\le n$, then there
  exists a unique $V\in L_2(D;\wedge^{k-1})$ such that $\del V=0$, $\nu\lctr V|_\Sigma=0$
  and
  $$
    F= d V.
  $$
  If $k=1$, then such $V$ exists, but is unique only up to constant scalar functions.  
\item[{\rm (ii)}]
  If $F\in L_2(D;\wedge^k)$ satisfies $\del F=0$, and $0\le k\le n-2$, then there
  exists a unique $V\in L_2(D;\wedge^{k+1})$ such that $dV=0$, $\nu\wedg V|_\Sigma=0$
  and
  $$
    F= \del V.
  $$
  If $k=n-1$, then such $V$ exists, but is unique only up to constant $n$-vector fields.
\end{itemize}
\end{prop}
In terms of operators, this result means that
$$
  d: \nul(\udel;\wedge^{k-1})\longrightarrow \nul(d;\wedge^k),\qquad \text{and}
  \qquad \del: \nul(\ud;\wedge^{k+1})\longrightarrow \nul(\del;\wedge^k),
$$
are injective (except when $k=1, n-1$ respectively) and surjective unbounded operators.
This follows from \cite[Theorem 1.3]{AMc} together with Poincar\'e's lemma.
Here $\nul(T;\wedge^k):= \nul(T)\cap L_2(D;\wedge^k)$.

\begin{defn}
  Given a two-sided harmonic field $U:D^\pm\rightarrow \wedge^k\R^n$, i.e. $\del dU=0= d\del U$,
  we say that $V_1: D\rightarrow \wedge^{k-2}\R^n$ and $V_2:D\rightarrow \wedge^{k+2}\R^n$ 
  are {\em Hodge type harmonic conjugates} to $U$ if
  $(d+\del)(V_1+U+V_2)=0$ and if $\nu\lctr V_1|_\Sigma=0$ and $\nu\wedg V_2|_\Sigma=0$.
\end{defn}

Our main result in this section is the following theorem, which shows that the boundary 
condition imposed on Hodge type harmonic conjugates yields a well posed problem.

\begin{thm}   \label{thm:hodgeconj}
  Assume that $D\subset \R^n$ is Lipschitz diffeomorphic to the unit ball, let $0\le k\le n$
  and let $U\in L_2(D;\wedge^k)$ be such that $dU, \del U\in L_2(D;\wedge)$
  and $\del d U=0= d\del U$.
  If $k=n-1$, assume that $dU=0$, and if $k=1$ assume that $\del U=0$.
  Then there exists Hodge type harmonic conjugates $V_1,V_2$ to $U$ in $D$ such that
  $\|V_1\|_{L_2(D)}\lesssim\|\del U\|_{L_2(D)}$ and $\|V_2\|_{L_2(D)}\lesssim\|dU\|_{L_2(D)}$.
  The conjugates are unique, except if $k=2$, when $V_1$ is unique modulo constants,
  and if $k=n-2$, when $V_2$ is unique modulo constants.
\end{thm}

\begin{proof}
  Apply Proposition~\ref{prop:hodge}(i) to $F= \del U$, and Proposition~\ref{prop:hodge}(ii) to $F=d U$.
\end{proof}

\begin{rem} 
(1)  
  In the complex plane, when $n=2$ and $k=0$, harmonic conjugates in general are unique
  modulo constants. In particular Hodge type conjugates coincide with Cauchy type conjugates, modulo constants.
  Note that all $V=V_2:D\rightarrow \wedge^2\R^2$ are normal on the boundary, since $\nu\wedg V|_\Sigma=0$.
  
(2)
  For general domains $D$, Cauchy type harmonic conjugates and Hodge type conjugates
  will not coincide in general when $n\ge 3$, not even when $k=0$. To see this, note that there is no reason
  for the Cauchy type conjugate
$$
  V(x)= \int_\Sigma (E(y-x)\wedg \nu(y)) h(y) d\sigma(y),\qquad h:\Sigma\longrightarrow \wedge^0\R^n
$$
  to satisfy $\nu\wedg V|_\Sigma=0$.
  However, they do coincide, for all $0\le k\le n$, when $D$ is a sphere. 
  To see this for the unit sphere, note that the normal vector $\nu(y)$ in this case is $y$,
  so
$$
  V_2(x)=\int_\Sigma E(y-x)\wedg \nu(y)\wedg h(y) \, d\sigma(y)= -x\wedg \int_\Sigma \frac{y\wedg h(y)\, d\sigma(y)}{\sigma_{n-1}|y-x|^n}
$$
 is normal on $\Sigma$. Similarly, the Cauchy type conjugate $V_1$ is seen to be tangential on
 $\Sigma$ in this case. Hence they coincide with the Hodge conjugates by uniqueness in Theorem~\ref{thm:hodgeconj}
 (modulo constants when $k=2,n-2$).
\end{rem}

\bibliographystyle{acm}
%GATHER{AKMcDirac.bib}  % makes sure WinEdt finds citations...
%\bibliography{hilbert}

\end{document}